\documentclass[12pt]{article}
\input epsf
\usepackage{epsfig}
\def\C{\mathbf{C}}
\def\R{\mathbf{R}}
\def\bC{\mathbf{\overline{C}}}
\def\K{K_{\mathrm{crit}}}
\def\k{k_{\mathrm{crit}}}
\def\Ima{\mathrm{Im}\, }
\def\Rea{\mathrm{Re}\, }
\def\id{\mathrm{id}}
\def\T{\mathbf{T}}
\def\Z{\mathbf{Z}}
\def\N{\mathbf{N}}
\def\S{\mathbf{S}}

\begin{document}
\title{On metrics of curvature $1$ with four conic singularities
on tori and on the sphere}
\author{Alexandre Eremenko\thanks{Supported by NSF grant DMS-1361836.}
\
and Andrei Gabrielov\thanks{Supported by NSF grant DMS-1161629.}}
\maketitle
\begin{abstract}
We discuss conformal metrics of curvature $1$ on tori and on the sphere,
with four conic singularities whose angles are multiples of $\pi$.
Besides some general results we study in detail the family of
such symmetric metrics on the sphere,
with angles $(\pi,3\pi,\pi,3\pi)$. As a consequence we find
new Heun's equations whose general solution is algebraic.


\end{abstract}

\section{Introduction} The general problem studied here concerns conformal
Riemannian metrics of constant positive curvature with conic singularities
on a compact Riemann surface $S$. For the history of the problem and its relevance
to mathematics and physics
we refer to \cite{Troy2,LT,Tarantello,Lin1,FKKRUY,EGT1,EGT2,MP}.
The goal is to classify such metrics up to isometry, in particular
{\em to determine for a given Riemann surface,
how many such metrics exist with prescribed singularities
and prescribed
angles at each singularity.}

Analytically the problem consists in the study of the Liouville
equation
\begin{equation}\label{liouville}
\Delta u+e^{2u}=0
\end{equation}
on the punctured Riemann surface, with prescribed singularities
at the punctures. If $\rho(z)|dz|$ is the length element of the metric
in a local conformal coordinate, then $\rho=e^u$, and equation (\ref{liouville})
expresses the fact that the curvature equals $1$. The behavior at a singularity
$a\in S$ is $\rho(z)\sim c|z|^{\alpha-1}$, where $c>0$,
$z$ is the local coordinate
which is $0$ at $a$, and $2\alpha>0$ is the angle at the singularity.

Here and in what follows we measure
the angles at the singularities in {\em half-turns},
for example ``integer angle'' means an integer multiple of $\pi$ radians.

An important special case is that of {\em symmetric} metrics, that is
metrics which are invariant with respect to an anti-conformal involution
which leaves all singularities fixed. In the case of the sphere,
such metrics are in one-to-one correspondence with spherical polygons
\cite{EGT2}.
So the problem in this case consists in classification and enumerating
spherical polygons of prescribed conformal type with prescribed angles
at the corners.

Every surface of constant curvature $1$ is locally isometric to a piece
of the standard sphere. This isometry is analytic in conformal coordinates,
so it admits an analytic continuation along every path not passing through
the singularities, and we obtain a multi-valued analytic function
$$f:S\backslash\{ a_0,\ldots,a_{n-1}\}\to\bC,$$
where $a_j$ are the singularities, and $\bC$ is the Riemann sphere.
This function $f$ is called the {\em developing map}, and the length element
of the metric is expressed in the form
$$\rho(z)|dz|=\frac{2|f'||dz|}{1+|f|^2}.$$
Monodromy of the developing map consists of linear-fractional transformations
(rotations of the sphere), so the Schwarzian derivative
$$\S[f]:=\frac{f'''}{f'}-\frac{3}{2}\left(\frac{f''}{f'}\right)^2$$
is single valued on $S\backslash\{ a_0,\ldots,a_{n-1}\}$.
The behavior of $f$ near the conic singularities $a_j$
implies that $\S[f]$
has double poles at these points, with principal parts
related to the angles (see below). This implies that
$\S[f]=R$, where $R$ is a rational function in the case
of the sphere, and an elliptic function in the case of a torus.
General solution of this Schwarz equation is a ratio
of two linearly independent solutions of the linear differential
equation with regular singular points
\begin{equation}\label{ode}
w''+\frac{1}{2}Rw=0.
\end{equation}
The exponent differences at the singular
points are the angles. After a normalization $n-3$ free parameters
remain which are called {\em accessory parameters}.
In terms of this equation, our general problem is equivalent to
the following:

{\em For given singular points and the exponent differences, how many
values of accessory parameters exist, so that the monodromy group
of (\ref{ode}) is conjugate to a subgroup of $PSU(2)$.}

In the case of three singularities, there is no accessory parameter,
the equation (\ref{ode}) is hypergeometric, and the problem of classification
of the metrics was completely
solved for this case in \cite{Klein,E,FKKRUY}. The case of four singularities
corresponds to Heun's equation which has one accessory parameter. The assumption
that the angles are integers leads to the special case of the
Heun equation, which is equivalent to the so-called
Darboux--Treibich--Verdier equation \cite{Veselov},
and which is easier
to study because its general solution is meromorphic on the universal cover
of the torus, and can be expressed
in terms of Abelian integrals.

For $n\geq 4$, the existing results show that the general
problem described above
is very difficult,
see, for example, the recent study of metrics with one singularity
of integer angle on a torus \cite{Lin1,Lin2,Lin3}, or
classification of such symmetric metrics with four angles, at least
one of them even integer, on the sphere \cite{EGT1,EGT2,EGT3}.

In this paper we study new special cases. We restrict ourselves
to the cases of four singularities with even angles on tori
and integer angles on the sphere.

In the first part of this paper we establish a relation between the metrics on
the sphere and on tori, and restate the results obtained in
\cite{EGSV,EGT1,EGT2,EGT3}. In the second part we study the simplest
family which is not covered by the existing results,
the metrics with angles $(6,2,6,2)$ on rectangular tori, satisfying
some additional symmetry conditions. This metric has in fact
only two conic
singularities. By the results
of the first part, this is equivalent to
the study of
symmetric metrics on the sphere with angles $(3,1,3,1)$.

We thank Walter Bergweiler, Chang-Shou Lin and Dmitrii Novikov
for useful discussions
of the problems considered here. We also thank the anonymous 
referee who corrected several
mistakes and suggested improvement of the exposition.

\section{Metrics on the sphere and on tori}
We consider a metric of curvature $1$ on a torus
$$\T=\C/\Lambda,\quad \Lambda=\{ m\omega_1+n\omega_3:m,n\in\Z\},\quad \tau=\omega_3/\omega_1,\quad\Ima\tau>0,$$
with four conic singularities
at half-periods $\omega_j/2,\; 0\leq j\leq 3$, where $\omega_0=0,\;
\omega_2=\omega_1+\omega_3,$
and the metric has an even angle at each singularity.
Let the angles be $2m_j,\; m_j\in \N$. Then the developing map $F$
must satisfy the Schwarz equation $\S[F]=R$, where $R$ is an elliptic
function of the form
$$R(z)/2=
-\sum_{j=0}^3\left(n_j(n_j+1)\wp(z-\omega_j/2)+b_j\zeta(z-\omega_j/2)\right)
+\lambda,
$$
with $n_j=(m_j-1)/2,$ and $\sum_{j=0}^3b_j=0$, as is explained below.
The general solution of the Schwarz equation is a ratio of two
linearly independent solutions of the equation
\begin{equation}\label{1}
w''-\left(\sum_{j=0}^3n_j(n_j+1)\wp(z-\omega_j/2)+\sum_{j=0}^3
b_j\zeta(z-\omega_j/2)+\lambda\right)w=0.
\end{equation}
The exponents at the singularities of (\ref{1})
are determined from the characteristic equation
$$\rho(\rho-1)-n_j(n_j+1)=0,$$
whose solutions are $\rho=1/2\pm(n_j+1/2)$, so the difference of the
exponents is
$m_j=2n_j+1$, and these are halves of
our prescribed angles of the metric on $\T$.

Monodromy of the developing map cannot contain a parabolic transformation,
so there should be no logarithms in the local solutions of (\ref{1}) near
the singularities, and equation
(\ref{1}) must have trivial monodromy at $\omega_j/2$.
In other words, all singularities must be {\em apparent}.
This gives three polynomial conditions on the parameters $b_j$, see,
for example \cite{GV}. In addition to this, we have conditions
on the accessory parameter $\lambda$ which ensure that monodromy
of (\ref{1}) is unitarizable, that is conjugate to
a subgroup of $SU(2)$. In particular, these conditions of
absence of logarithms are always satisfied if all $b_j=0$, in which case
(\ref{1}) becomes the Darboux--Treibich--Verdier equation
\begin{equation}\label{dtv}
w''-\left(\sum_{j=0}^3n_j(n_j+1)\wp(z-\omega_j/2)+
\lambda\right)w=0.
\end{equation}
For this equation
an explicit formula for the general solution of
(\ref{1}) can be written \cite{VV}.

The monodromy coming from the fundamental
group of the torus must be conjugate to
an Abelian subgroup of $PSU(2)$. There are
only two kinds of such groups \cite{Klein2}:
the Klein Viergroup $(\Z/2\Z)^2$ and a cyclic group
(finite or infinite). By a post-composition of $F$ with an element
of $PSU(2)$, which does not change the metric,
we can achieve that $F$ satisfies
one of the following:
\begin{equation}\label{A}
F(z+\omega_1)=-F(z),\quad F(z+\omega_3)=1/F(z),
\end{equation}
or
\begin{equation}\label{B}
F(z+\omega_j)=e^{2\pi i\alpha_j}F(z),\quad j=1,3.
\end{equation}
We call two metrics {\em equivalent}
if their developing maps are related by post-composition with
a linear-fractional transformation.
Thus in the case (\ref{B}) each metric is a member of a
one-parametric
family of equivalent metrics with developing maps $kF$, $k\neq 0$.
\vspace{.1in}

\noindent
{\bf Theorem 1.} {\em Suppose that $b_j=0,\; 0\leq j\leq 3$ in
(\ref{1}).
Then in case (\ref{A}) the metric is even,
that is invariant under $z\mapsto-z$,
while in case (\ref{B}) each
equivalence class contains an even metric.}
\vspace{.1in}

This generalizes results from \cite{Lin2}, sections 2.1 and 5.1.
\vspace{.1in}

{\em Proof.} It is sufficient to show that
\begin{equation}\label{2}
F(-z)=TF(z),
\end{equation}
with $T\in PSU(2)$.

That $F$ satisfies (\ref{2}) with {\em some} $T\in PSL(2,\C)$ is clear because
the Schwarzian $\S[F]$ is even.
It is also clear from (\ref{2}) that
\begin{equation}\label{idem}
T^2=\id,
\end{equation}
so $T$ is either an elliptic transformation with multiplier $-1$ or the identity.
If $T=\id$, the proof is finished, so assume that the multiplier of $T$
is $-1$.

Suppose now that $F$ is a meromorphic function, satisfying (\ref{2}), and
\begin{equation}\label{S}
F(z+\omega)=SF(z)
\end{equation}
where $S$ is a linear-fractional transformation. We claim that in this case
the equation
\begin{equation}\label{com}
STS=T
\end{equation}
must hold.
Indeed, (\ref{S}) implies
\begin{equation}\label{inS}
F(z-\omega)=S^{-1}F(z)
\end{equation}
and also
\begin{equation}\label{SS}
F(-z+\omega)=SF(-z).
\end{equation}
Using (\ref{2}), (\ref{SS}) and (\ref{inS}) we obtain
$$STF(z)=SF(-z)=F(-z+\omega)=F(-(z-\omega))=TF(z-\omega)=TS^{-1}F(z),$$
which implies (\ref{com}).

Suppose that (\ref{A}) holds. Both equations in (\ref{A}) are of
the form (\ref{S}) with $S^{-1}=S$, thus $ST=TS$. Using the first equation
(\ref{A}) in which $S(z)=-z$, we conclude that $T(z)=cz$ or $T(z)=c/z$,
$c\neq 0$. In the first case (\ref{idem}) implies $c^2=1$, so $T\in PSU(2)$.
In the second case, we use the second equation in (\ref{A}) which implies
that $T$ commutes with $S(z)=1/z$, and conclude again that $c^2=1$ and $T\in PSU(2)$.

Suppose now that (\ref{B}) holds, and $S(z)=kz,\; k=e^{2\pi i\alpha},\; k\neq 1$.
Then (\ref{com}) gives
\begin{equation}\label{kwa}
kT(kz)\equiv T(z).
\end{equation}
Substituting
$$T(z)=\frac{az+b}{cz+d}$$
to (\ref{kwa}), after simple computation we obtain
$$bd=0\quad\mbox{and}\quad ac=0,$$
and
$$ad=0\quad\mbox{or}\quad k=-1.$$
If $b=0$, then $a\neq 0,\; d\neq 0,$ so $k=-1$ and $T(z)=(a/d)z$. Using (\ref{idem}) we conclude that $T(z)=-z$.

If $d=0$, then $c\neq 0,\; b\neq 0$, so $a=0$ and $T(z)=b/(cz)$,
and after multiplication of $F$ by a constant we obtain $T\in PSU(2)$.

Finally, if we have (\ref{B}) with both multipliers $k_j=1$, then $F$ is an elliptic function,
with $F(-z)=TF(z)$, where $T$ is an elliptic
transformation satisfying (\ref{2}). Let
$V\in PSL(2,\C)$ be such that $VTV^{-1}\in PSU(2)$
then the metric with the developing map $VF$
is even.

This completes the proof of Theorem 1.
\vspace{.1in}

Every even metric on the torus is a pull-back
of a metric from the sphere via the $\wp$ function.
Thus we obtain a metric on the sphere,
which has $4$ singularities at the critical
values of $\wp$, that is at $\infty, e_1,e_2,e_3$,
and the angles at these singularities are $2\alpha_j=m_j$.
We denote by $f$ the developing map on the sphere, so that $F=f\circ\wp$.
The developing map $f$ is a ratio of two solutions of the Heun equation
\begin{equation}\label{heun}
w''+\frac{1}{2}\left(\sum_{j=0}^2\frac{1-2n_j}{z-a_j}\right)w'+
\frac{\alpha'\alpha''z-E}{\prod_{j=0}^2(z-a_j)}w=0,
\end{equation}
where
$$\alpha'=(1+n_3-n_2-n_1-n_0)/2,$$
$$\alpha''=-(n_3+n_2+n_1+n_0)/2,$$
and $(a_0,a_1,a_2)=(e_2,e_3,e_1).$
Now we consider four cases.
\vspace{.1in}

\noindent
{\bf Case 1.} {\em All $m_j$ are odd.}
\vspace{.1in}

Then the monodromy of $f$ is generated by four non-trivial involutions.
\vspace{.1in}

\noindent
{\bf Proposition 1.} {\em In Case 1, $F$ satisfies (\ref{B}).}
\vspace{.1in}

This generalizes a result \cite[Thm. 0.3]{Lin2}.

{\em Proof.} Proving this by contradiction, assume that we have (\ref{A}).
Then the monodromy group of $F$ is the Klein Viergroup $V= (\Z/2\Z)^2$.
Let $\sigma_j$ be conformal involutions of $\C$ with fixed points
$\omega_j/2$,
$$\sigma_j(z)=\omega_j-z,\quad 0\leq j\leq 3.$$ Their images in the monodromy
group are non-trivial elements of $V$,
because all $m_j$ are odd.
As $V$ contains only three non-trivial elements, two of
the distinct involutions $\sigma$ and $\sigma'$ have the same image,
so the image
of $\sigma\sigma'$ in the monodromy group must be the identity.
But $\sigma\sigma'$ is a translation by some $\omega_j$, so we obtain
a contradiction with (\ref{A}), which proves the Proposition.
\vspace{.1in}

In the remaining cases, some $m_j$ are even, thus the metric
on the sphere has at least one even angle. All these cases were studied in
detail and completely classified in \cite{Sch,EGSV,EGT1,EGT2,EGT3}.
We state the results.
\vspace{.1in}

\noindent
{\bf Case 2.} {\em Exactly one of the $m_j$ is even.}
\vspace{.1in}

Let $m_0=2\alpha_0$ be the even angle. Then it follows from \cite[Proposition 4.2]{EGT3} that
there are at most $\alpha_0$ distinct metrics, and for generic tori
there
are exactly $\alpha_0$ of them.
The accessory parameters $\lambda$ are found
from the algebraic equation that expresses the absence of logarithms in
the solution of the Heun equation.
\vspace{.1in}

\noindent
{\bf Case 3.} {\em Exactly two of the $m_j$ are even.}
\vspace{.1in}

Let $m_j=2\alpha_j$, $0\leq j\leq 3$ and assume that
$$\alpha_3=r+1/2,\quad\alpha_0=s+1/2,\quad \alpha_1=p,\quad\alpha_2=q.$$
Then according to \cite[Thm. 4.1]{EGT1} there are two
subcases:

Subcase 3a. $p+q+r+s$ is even, $|r-s|\leq p+q$,
and the number of classes of metrics
is at most $(p+q-|r-s|)/2$ with equality for generic tori,

Subcase 3b. $p+q+r+s$ is odd, $r+s+3\leq p+q$, and the number of classes
of metrics is at most $(p+q-r-s-1)/2$ with equality for generic tori.

If three of the $m_j$ are even then all four must be even, because the
product of local monodromies is the identity,
and we obtain the
last case:
\vspace{.1in}

\noindent
{\bf Case 4.} {\em All $m_j$ are even.}
\vspace{.1in}

In this case, $f$ is a rational function with four critical points of orders
$m_j/2$. Such functions were classified in \cite{Sch,EGSV}.

So in all cases 2,3,4, we have a complete classification
of metrics in question.

In the next sections we consider a special family which belongs to Case~1.

\noindent
\section{Spherical rectangles. Statement of results}
Here we begin to study spherical
quadrilaterals whose angles are odd multiples of $1/2$. We call them
{\em spherical rectangles}.

Metrics on the sphere corresponding to spherical rectangles are symmetric
with respect to a circle, and without loss of generality we may assume
that this circle is the real line. Singularities lie on the real line.
Consider the ramified covering of degree $2$
from a torus to the sphere, ramified over
these singularities. It is a composition of an element of $PSL(2,\R)$
with the Weierstrass $\wp$ function whose invariants $e_j=\wp(\omega_j/2)$
are real. Such Weierstrass functions are defined on ``rectangular'' tori,
which means that $\omega_3/\omega_1=iK,\; K>0$.
The pull back of the metric on the torus has singularities with even
angles
at half-periods $\omega_j/2$.

In each quadrilateral we mark one corner. Two quadrilaterals
are considered congruent if there is an orientation preserving
isometry sending one to another
and sending the marked corner to the marked corner.

We usually choose the upper half-plane $H$ as a conformal model of a spherical
polygon. Then the corners $a_j$ become points on the real line,
and the metric has the length element $\rho(z)dz$, $\rho=e^u$ where $u$ satisfies
(\ref{liouville}). Then the spherical polygon is a pair $Q=(H,\rho)$,
where $\rho$ has conic singularities at $a_j$ and the sides
$(a_j,a_{j+1})$ are geodesic with respect to the metric $\rho$.

Developing map of a quadrilateral satisfies the
Schwarz equation associated with the Heun equation, as explained
in the Introduction (see also \cite{EGT1}).

Our first result is
\vspace{.1in}

\noindent
{\bf Theorem 2.} {\em Let $f:Q\to \bC$ be the developing map of a
spherical quadrilateral whose angles are odd multiples of $1/2$.
Then there are two opposite sides whose $f$-images are
contained in the same circle. On the other hand $f(\partial Q)$ is
not contained in two circles.}
\vspace{.1in}

Thus there are two types of such spherical quadrilaterals: in the first type
one of the two sides whose images belong to the same circle ends at the
marked corner, and in the second type such a side begins at the marked corner.

To state our main result we define the conformal modulus.
Let us map conformally our quadrilateral onto the rectangle
$(0,1,1+iK,iK)$, so that the marked corner is mapped to $0$. Then $K>0$ is
called the modulus.
\vspace{.1in}

\noindent
{\bf Theorem 3.} {\em There are two continuous families of
spherical quadrilaterals
with angles $(3/2,1/2,3/2,1/2)$, and angle $3/2$ at the marked corner.
One family consists of quadrilaterals
of the first type, and the modulus in this family
varies from $0$ to $\K<1$.
The second family consists of quadrilaterals
of the second type, and the modulus
varies from $1/\K$ to $\infty$. Each family contains exactly one
quadrilateral of the given modulus in the described range.}
\vspace{.1in}

So there is a ``forbidden interval'' $[\K,1/\K]$ for the moduli of
spherical quadrilaterals with angles $(3/2,1/2,3/2,1/2)$.
\vspace{.1in}

Theorem 3 is proved in the next section.

{\em Proof of Theorem 2.}
Let us enumerate the great circles $C_j$ which contain the images of sides
according to positive orientation
of the boundary, so that $j$ is a residue modulo $4$.
We may assume that $C_1$ and $C_2$ are real and imaginary axes,
then $C_3$ must cross the real axis perpendicularly at the points $r,-1/r$,
and $C_4$ must cross the imaginary axis perpendicularly at the points $it,-i/t$.
Without loss of generality $r\geq 0$ and $t\geq 0$.

If the $4$ circles are in general position (that is there are
no triple intersections),
then both $r$ and $t$ are finite, positive,
and there is a crossing point $A$
of $C_3$ and $C_4$. The centers of $C_3$ and $C_4$ are at $c_3=(r-1/r)/2$
and $c_4=i(t-1/t)/2$. Their radii are $r_3=(r+1/r)/2$ and $r_4=(t+1/t)/2$.
Then we obtain from the right triangle $(c_2,0,c_4)$ that
the square of the distance between the centers is
$$(r-1/r)^2/4+(t-1/t)^2/4.$$
On the other hand the angle between the radii of $C_3$ and $C_4$ at $A$
is right, so the same distance is
$$(r+1/r)^2/4+(t+1/t)^2/4,$$
a contradiction. Therefore the images of two opposite sides must lie on
the same circle.

The non-generic cases are easy and we leave them to the reader.

To prove the last statement of Theorem 2 by contradiction,
suppose that $f(\partial Q)$ is contained in two circles.
By post-composition of $f$ with a rotation of the sphere, we may assume
that these two circles are intersecting at $0$ and $\infty$,
and one of them is the real line.
As all angles of our quadrilateral are odd multiples of $1/2$,
the corners must be mapped to $0$ and $\infty$, and the other circle
is the imaginary line. Then $f^2$ extends by symmetry to a rational
function with two critical values $0$ and $\infty$. So we must have
$f^2(z)=cz^n$, but then $f^{-1}(\{0,\infty\})$ consists
of only two points. This is a contradiction which proves the last statement
of Theorem 2.
\vspace{.1in}

Let our three circles be the real line, the line $\ell_\alpha$ through
the origin making angle $\alpha$ with positive ray, and the unit circle.
The $f$-images of the two opposite sides lie on the
unit circle and the images of the other two
sides on $\R$ and $\ell_\alpha$, respectively. We have one continuous
parameter $\alpha\in(0,1)$, and the modulus of the quadrilateral
depends on the angles and on $\alpha$.

The group generated by reflections in the $f$-images of the
sides has a subgroup of index $2$
consisting of conformal maps; this subgroup is the monodromy
group of our developing map $H\to\bC$. The monodromy group
is generated by $3$ elliptic
involutions
with fixed points $(0,\infty)$, $(1,-1)$ and $(e^{\pi i\alpha},-e^{\pi i\alpha})$,
so it consists of the maps of the form
\begin{equation}\label{monodromy}
z\mapsto e^{2\pi i m\alpha}z\quad\mbox{and}\quad
z\mapsto e^{2\pi i m\alpha}/z,
\end{equation}
where $m$ is any integer.
When $\alpha$ is rational, this group is finite, it is a dihedral group,
and the developing map is algebraic.
For this case, there is a theorem of Klein (see \cite{BD})
which says that the Heun equation
is a pullback of a hypergeometric equation with finite monodromy group.
This means that the developing map factors as $h\circ b$, where $h$ is
a solution of the hypergeometric equation (in our case $h^{-1}$ is a conformal
map of the triangle with angles $(1/2,1/2,2\alpha)$ onto the upper half-plane,
and $b$ is a Belyi rational function which can be computed algebraically).


Our monodromy group (\ref{monodromy}) belongs to the class
of groups considered by Van Vleck which consist of transformations
of the form
\begin{equation}\label{vanvleck}
z\mapsto cz,\quad z\mapsto d/z.
\end{equation}
Van Vleck \cite{VV} proved that Fuchsian differential equations with such
projective monodromy groups can be characterized by the property that
they have two linearly independent solutions whose product
is a polynomial. This remarkable property was first discovered by Hermite
\cite{Hermite} for the case of Lam\'e equation, and generalized by Darboux
\cite{Darboux} to equation (\ref{dtv}) which is equivalent to (\ref{heun}).
It permits to solve
the equation in a closed form.
If $P$ is this polynomial, then two linearly independent solutions
of (\ref{heun}) can be represented in the form
\begin{eqnarray}\label{vv}
w_1&=&\sqrt{P}\exp\left( C\int \prod(z-a_j)^{\alpha_j-1}\frac{dz}{P(z)}\right),\\
w_2&=&\sqrt{P}\exp\left(-C\int \prod(z-a_j)^{\alpha_j-1}\frac{dz}{P(z)}\right).
\end{eqnarray}
Hermite found the method of finding $P$ in terms of $a_j,\alpha_j$ and the
accessory parameter.
Equation on the accessory parameter comes from the condition
that the monodromy group has the special form (\ref{monodromy}).
In section 4 we perform all computations explicitly in a simple
case.

To pull back our metric on a torus,
we model our
spherical quadrilateral on a flat rectangle, instead of the half-plane.
Consider the rectangle
$(0,\omega_1/2,\omega_2/2,\omega_3/2)$, where $\omega_1$ and $\omega_3/i$ are
positive and $\omega_2=\omega_1+\omega_3$.
Now the developing map $f$ is defined on this rectangle
and maps it to the standard sphere. Without loss of generality
one can assume that
$\omega_1=2$ and $\omega_3=2iK$, but sometimes we will
use $\omega_1$ and $\omega_3$. Then the developing map extends by reflection
to the whole
plane. We assume without loss of generality
that the vertical sides of the rectangle are mapped by $f$ to
the unit circle, and the side $(0,\omega_1/2)$
to the real line.
Then the developing map is real on the real line and
satisfies
\begin{equation}\label{periods}
f(z+\omega_1)=f(z),\quad f(z+\omega_3)=e^{2\pi i\alpha}f(z),
\end{equation}
which is a special case of (\ref{B}).
To this one can add
\begin{equation}\label{oddness}
f(-z)=1/f(z)
\end{equation}
and similar properties
for conformal involutions with fixed points at $\omega_j/2$.
If the angles $\alpha_j$ of our spherical rectangle are $n_j+1/2$,
then the map
$f$ is ramified at the corners and is $m_j=(2n_j+1)$-to-$1$ near $\omega_j/2$,
that is $f'$
has a zero of order $m_j-1$ at $\omega_j/2$ and no other zeros.

In (\ref{periods}), $\alpha$ and the parameter $\tau=\omega_3/\omega_1=iK$
are connected by some relation
which also depends on the $m_j$.

Thus we have the problem:
\vspace{.1in}

{\em For given periods $\omega_j$ and given integers $m_j$,
how many functions $f$ exist
(up to proportionality) that satisfy (\ref{periods})
and have prescribed critical point pattern: critical points of
order $m_j-1$ at $\omega_j/2$ and no other critical points.}
\vspace{.1in}

The Schwarzian derivative of $f$
is doubly periodic with periods $2\omega_1,2\omega_3$, and has double poles
at $\omega_j$, where we must have apparent singularities
and the exponent difference must be $m_j$. So our developing map
is the ratio of two solutions
of the Darboux--Treibich--Verdier equation (\ref{dtv}).
Thus we have an equivalent problem:
\vspace{.1in}

{\em For given periods $\omega_j$ and given integers $n_j$, how many
real values of $\lambda$
exist, so that a ratio $f=w_1/w_2$ of some linearly independent solutions
of (\ref{1}) has properties (\ref{periods}) with some real $\alpha$.}
\vspace{.1in}

\noindent
Parameter $\lambda$ must be real because
the Schwarzian of $f$ is real on the real line, and
$\wp(z-\omega_j/2)$ are real on the real line for the lattices with
real $\omega_1$ and imaginary $\omega_3$.
\vspace{.1in}

\noindent
\section{Spherical rectangles with angles $(\frac32,\frac12,\frac32,\frac12)$}
\vspace{.1in}

In this section we prove Theorem 3.
We consider the simplest spherical rectangles, with the angles $(3/2,1/2,3/2,1/2)$.
The angles are measured in half-turns, so that $1/2$ means $\pi/2$ radians.
\vspace{.1in}

\noindent
{\bf Lemma 2.} {\em There are two families of spherical rectangles
with angles
\newline
$(3/2,1/2,3/2,1/2)$, each depending on a parameter
$\alpha\in(0,1)$. Rectangles of these families
are isometric to subsets of the sphere whose oriented boundaries
are shown in Fig.~1.
}
\vspace{.1in}

{\em Proof.}
The area of a rectangle with angles as in Lemma 2 is
$1/2$ of the area of the sphere; this follows from the Gauss--Bonnet
formula.
We can normalize the developing map so that the image
of the boundary is contained in the union of the real line,
the line $\ell_\alpha$ as in the previous section, and
the unit circle.

The image of any side cannot cover more than $1/2$ of the great circle
to which this image belongs, otherwise this image would
contain either a half-plane or the exterior or the interior of the unit disk.
Since the area of the whole rectangle is $1/2$ of the area of the sphere,
it would follow that the rectangle is isometric to a hemisphere, which is
of course impossible.

From Theorem 1 we know that some pair of opposite sides is mapped by
the developing map to the same circle. Without loss of generality we assume
that this circle is the unit circle.
Let $s,s'$ be these two sides.
They are oriented in the usual way,
so that the rectangle is on the left of each side.
The angle at the beginning
of $s$ is the same as that at the beginning of $s'$,
and it is either $1/2$ or $3/2$.
So we have (at least) two types of such rectangles.

Consider a rectangle of the first type.
Let the sides be $(d,s,d',s')$ in the
order of positive orientation. Without loss of generality, the image of $d$
is on the real line, oriented in the increasing direction. This image must be either
$(-1,1)$ or its complement. Without loss of generality, let it be the complement.

Then the image
of $s$ must be a clockwise arc of the unit circle from $-1$ to some point $e^{\pi i\alpha}$
in the upper half-plane. Tracing the rest of the boundary we arrive
at a picture in the top of Fig.~1.

Rectangles of type 2 are considered in the similar way.
This proves the lemma.

An alternative way to understand spherical quadrilaterals is through
their {\em nets} \cite{EG,EGT2,EGT3}.
The developing map $f$ maps the sides of a quadrilateral $Q$
to three transversal great circles on the sphere.
The {\em net} $\Gamma$ of $Q$ is the preimage in $Q$ of these three circles.
For each circle $C$, its preimage in $Q$ is called $C$-{\em net}, denoted $\Gamma_C$.
An {\em arc} of the net $\Gamma_C$ is a connected component
of $\Gamma_C\setminus\partial Q$.
The net $\Gamma$ defines a triangulation of $Q$, each face of it being
mapped by $f$ one-to-one onto one of the triangles into which the three
circles partition the sphere.
\vspace{.1in}

\noindent{\bf Lemma 3.} {\em Let $C$ be the circle to which two opposite
sides $s$ and $s'$ of a quadrilateral $Q$ with angles $(3/2,1/2,3/2,1/2)$
are mapped. Then there is an arc $\gamma$ of $\Gamma_C$ with
the ends at the two corners
of $Q$ having angles $3/2$.}
\vspace{.1in}

{\em Proof.}
Since the area of $Q$ is $1/2$ of the area of the sphere,
and each corner having angle $3/2$ has three faces of $\Gamma$ adjacent to it,
with the total area greater than $1/4$ of the area of the sphere,
there should be a face $F$ of $\Gamma$ adjacent to both of these corners.
Thus these corners are connected by one of the sides of $F$ which must
be an arc of $\Gamma_C$ because $C$ is the only circle to which the
images of both of these corners belong.
\vspace{.1in}

\noindent{\bf Corollary.} {\em Each quadrilateral $Q$ with angles
$(3/2,1/2,3/2,1/2)$
is a union of two isometric halves of hemispheres,
each of them bounded by the arcs of two great circles
intersecting at a right angle.}

\vspace{.1in}
\noindent{\bf Remark.}
The statement of Lemma 3 is a special case of the following statement,
which will be proved elsewhere:
\vspace{.1in}

{\em Let $Q$ be a spherical $n$-gon such that
all its sides are mapped to three transversal great circles by the
developing map, and all its corners are mapped to some intersection points
of those circles. Then there is (in general, non-unique) triangulation of $Q$ by $n-3$ disjoint arcs of its net,
each arc having both ends at the corners of $Q$.}

\vspace{.1in}
We recall that our quadrilateral has a
marked corner, and two quadrilaterals are considered congruent if there is
an isometry between them which preserves the orientation and sends
the marked corner to the marked corner.

We choose the marked corner so that the angle at it is $3/2$,
and the two families are distinguished by the property
that the side ending at the marked
corner in the first family has the image on the same circle as its
opposite side, while in the second family the side that ends
at the marked corner and its opposite side have
images on two different circles.

Note that the choice of one of the two corners with the angle $3/2$ as marked
is not important, as our families have isometries $z\mapsto e^{\pi i\alpha}/z$
and $z\mapsto e^{-\pi i\alpha}/z$, respectively, exchanging the two corners.

Let us map our quadrilateral onto a flat rectangle with two sides $[0,1]$ and $[0,Ki]$,
where $0$ is the marked corner. Then $K>0$ is the {\em modulus} of the
quadrilateral.

\begin{figure}
\centering
\includegraphics[width=4in]{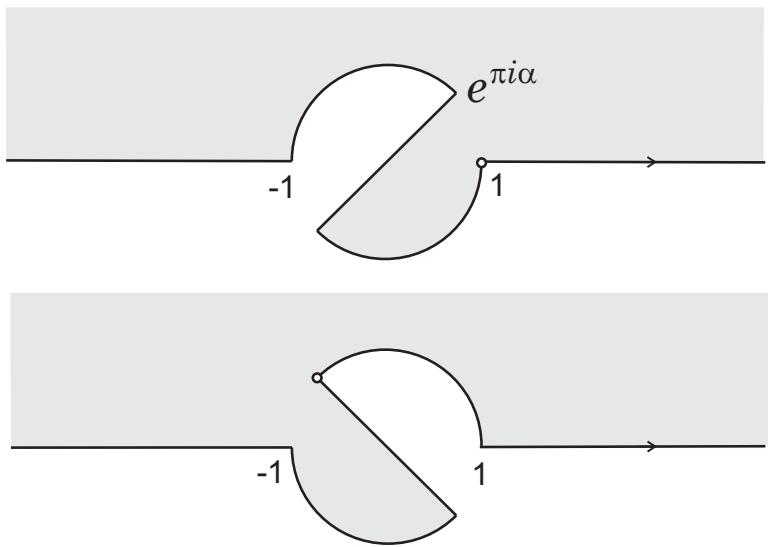}
\caption{Two families of quadrilaterals}\label{2families}
\end{figure}

\begin{figure}
\centering
\includegraphics[width=4in]{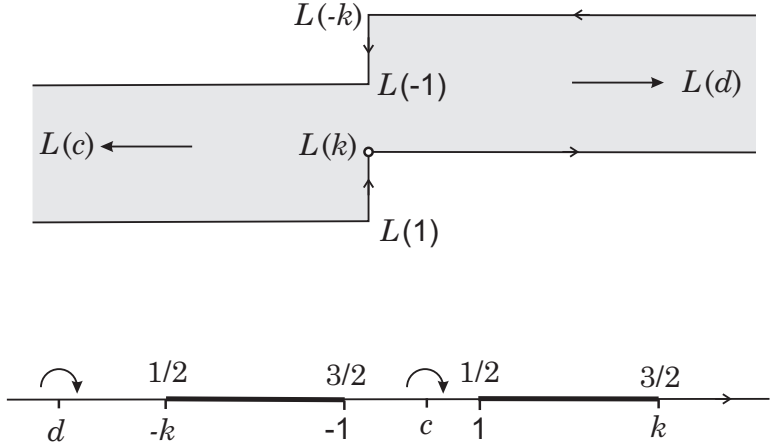}
\caption{Conformal map, $L=\log f$.}\label{rect}
\end{figure}

Now we return to the proof of Theorem 3.
Our two families are shown in Fig.~\ref{2families}. The marked corner is circled.
Consider the first family. Taking log we obtain the region in Fig.~\ref{rect}.
Now according to Schwarz and Christoffel, logarithm of the developing map
defined in the
upper half-plane is
$$L(z):=\log f(z)=A\int_k^z\sqrt{\frac{(1+\zeta)(k-\zeta)}{(1-\zeta)(k+\zeta)}}
\frac{d\zeta}{(\zeta-c)(\zeta-d)}.$$
This should be compared with Van Vleck's general formula (\ref{vv}).
The residue at $c$ equals:
$$\frac{A}{c-d}\sqrt{\frac{(1+c)(k-c)}{(1-c)(k+c)}}=\pm1.$$
Then the residues at $d$ are $\mp1$ by the Residue theorem. This gives
the condition
\begin{equation}\label{g}
h(c)=h(d)\quad\mbox{where}\quad h(x)=\frac{(1+x)(k-x)}{(1-x)(k+x)}.
\end{equation}
This is the ``Bethe-Ansatz equation'' for this case.
Thus $d=-k/c$ and
\begin{equation}\label{AA}
A=(c+k/c)\sqrt{\frac{(1-c)(k+c)}{(1+c)(k-c)}}>0.
\end{equation}
Therefore our integral is
\begin{equation}\label{3}
L(z)=\int_k^z\frac{c+k/c}{\zeta+k/c}
\sqrt{\frac{(1-c)(k+c)(1+\zeta)(k-\zeta)}{(1+c)(k-c)(1-\zeta)(k+\zeta)}}
\frac{d\zeta}{\zeta-c}=:\int_k^z\frac{g(c,\zeta)d\zeta}{\zeta-c},
\end{equation}
where the branch of the square root in the upper half-plane is chosen
so that the integrand is positive on $(k,\infty)$.
This and other similar integrals are analytic for $z$ in the upper half-plane,
and when $z$ is real we understand them as limits values of
this analytic function in the upper half-plane.

Here $c$ plays the role of the accessory parameter, and is simply related
to it.
The condition is that
$$\Rea \int_{-1}^k\frac{g(c,\zeta)d\zeta}{\zeta-c}=0,$$
which expresses the property that the points $L(k)$ and $L(-1)$
are on the same vertical line in Fig.~\ref{rect}. As the part of the integral (\ref{3})
over $[-1,1]$ is real and the part over $[1,k]$ is imaginary,
our condition is equivalent to
$$\int_{-1}^1\frac{g(c,\zeta)d\zeta}{\zeta-c}=0,$$
because the part of the integral in (\ref{3}) over $[-1,1]$ is real,
and its part over $[1,k]$ is imaginary.
{}From this condition $c\in(-1,1)$ must be found.
To regularize our integral on $(-1,1)$.
write it as
$$F(k,c):=\int_{-1}^1 (g(c,\zeta)-g(c,c))\frac{d\zeta}{\zeta-c}+
g(c,c)\int_{-1}^1\frac{d\zeta}{\zeta-c}.$$
Now $g(c,c)=1$ and
$$\int_{-1}^1\frac{dx}{x-c}=\log\frac{1-c}{1+c}.$$
The equation for $c$ is $F(k,c)=0$.

We would like to know that $c\in(0,1)$ is uniquely determined by $k$.
This would mean that there is at most one quadrilateral with prescribed
modulus in each family.
\vspace{.1in}

{\bf Proposition 1.} {\em There exists $\K\in(0,1)$ with the following property.
For every $K\in(0,\K)$ there exists unique $c\in (0,1)$ such that
$F(k,c)=0$. Moreover, $c\to 0$ as $K\to 0$ and $c\to 1$ as $K\to \K$.}
\vspace{.1in}

Here $K$ is a known (increasing) function of $k$ \cite{Ak}
see also (\ref{K}) below.
\vspace{.1in}

{\em Proof of Proposition 1}.
To prove uniqueness of $c$ it is sufficient to show
that after multiplication of $F$ by some function of $k$ and $c$ which is
non-zero in the considered range, the result is strictly decreasing in $c$
for each $k$.
First we consider the
function
\begin{eqnarray*}
F_1(k,c)&=&\sqrt{\frac{(1+c)(k-c)}{(1-c)(k+c)}}F(k,c)\\
&=&\int_{-1}^1\sqrt{\frac{(1+z)(k-z)}{(1-z)(k+z)}}\frac{(c+k/c)dz}{(z-c)(z+k/c)}.
\end{eqnarray*}
Then we do partial fraction decomposition
$$\frac{c+k/c}{(z-c)(z+k/c)}=
\left(\frac{1}{z-c}-\frac{1}{z+k/c}\right).$$
The second partial fraction $-1/(z+k/c)$ is strictly decreasing in $c\in(0,1)$, for each $z\in(-1,1)$, so it remains to prove that
$$F_2(k,c)=\int_{-1}^1\sqrt{\frac{(1+z)(k-z)}{(1-z)(k+z)}}\frac{dz}{z-c}$$
is decreasing.

\begin{figure}
\centering
\includegraphics[width=4in]{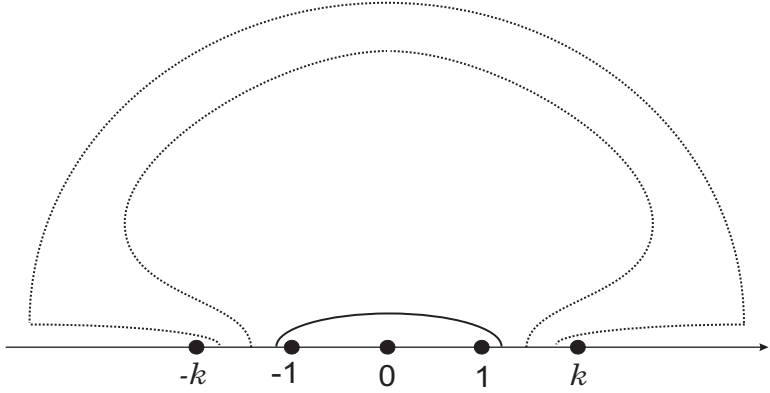}
\caption{Deformation of the path of integration.}\label{path}
\end{figure}

To do this we deform the contour of integration (see Fig.~\ref{path}). We consider the simple closed
curve $\gamma$ on the Riemann surface of the integrand which goes
once around the
segment $[-1,1]$ in the negative direction (clockwise). Then
$$F_2(k,c)=\frac{1}{2}\int_\gamma\sqrt{\frac{(1+z)(k-z)}{(1-z)(k+z)}}\frac{dz}{z-c}.$$
Then we deform the contour so that the part in the upper half-plane consists
of the interval $(k,R)$, $R>k$ oriented from $R$ to $k$, the interval
$(-R,-k)$ oriented from $-k$ to $-R$ and a large half-circle $|z|=R,\Ima z>0$.
The integral over the half-circle is $\pm\pi i$ times residue at infinity,
thus it is pure imaginary.

On the interval $(k,R)$ the integral has the form
$$-\int_k^R \phi(k,x)\frac{dx}{x-c},$$
where $\phi>0$ and it is decreasing as a function of $c$.
On the interval $(-R,-k)$ the integral has similar form, and also decreasing
as a function of $c$. This proves that for every $k$ there exists
at most one $c$ such that $F(k,c)=0$. Let this value be $c(k)$.

When $k\to 1+$, we have $c(k)\to 0+$ and $c(k)$ is increasing. So we have
some interval $(0,\k)$ on which $c$ is increasing.

When $k$ passes $1$, $c$ passes $0$, and we obtain the second family.

\begin{figure}
\centering
\includegraphics[width=4in]{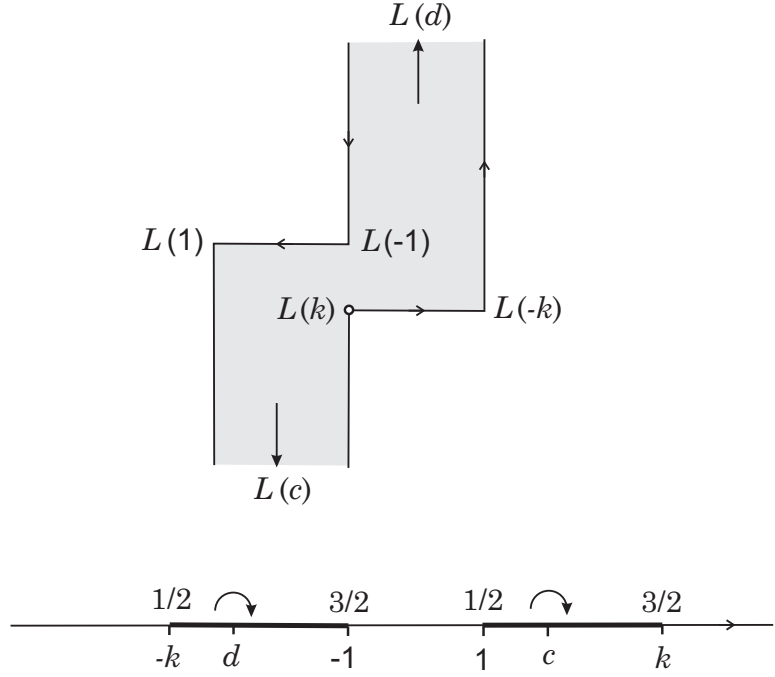}
\caption{Conformal map $L=\log f$ for $k>k_{\mathrm{crit}}$.}\label{non-geo}
\end{figure}

When $k>k_{\mathrm{crit}}$, the quadrilateral 
must have the shape as in Fig.~\ref{non-geo}.
The conformal map is
$$L(z)=\log f(z)=A'\int_k^z\sqrt{\frac{(1+\zeta)(k-\zeta)}{(1-\zeta)(k+\zeta)}}
\frac{d\zeta}{(\zeta-c)(\zeta+k/c)}=:A'\int_k^z G(\zeta,c)d\zeta.$$
We want to keep $A'>0$ but the residue at $c$ now must be $i$, and the
residue at $d$ must be $-i$.
That is
$$\frac{A'}{c-d}\sqrt{\frac{(1+c)(k-c)}{(1-c)(k+c)}}=\pm i.$$
Notice that the equation for $d$ is exactly the same as before,
thus $d=-k/c$, and that the expression under the root is now negative, and
we obtain
$$A'=(c-d)\sqrt{\frac{(c-1)(k+c)}{(c+1)(k-c)}}>0,$$
as desired.
The condition now is that
$$\Re A\int_{-1}^kG(\zeta,c)d\zeta=0\;\leftrightarrow \int_{-1}^1G(\zeta,c)d\zeta=-\pi,$$
and this must define $c\in(1,k)$.
The integral is now
$$\int_{-1}^1\frac{c^2+k}{cx+k}\sqrt\frac{(c-1)(k+c)(1+x)(k-x)}{(c+1)(k-c)(1-x)(k+x)}\frac{dx}{x-c}.$$
The only difference in comparison with the previous integral is
that factor $(1-c)$ is replaced by $(c-1)$ under the square root,
and the integral must be equal to $-\pi$. This integral does not require
desingularization.

\begin{figure}
\centering
\includegraphics[width=4in]{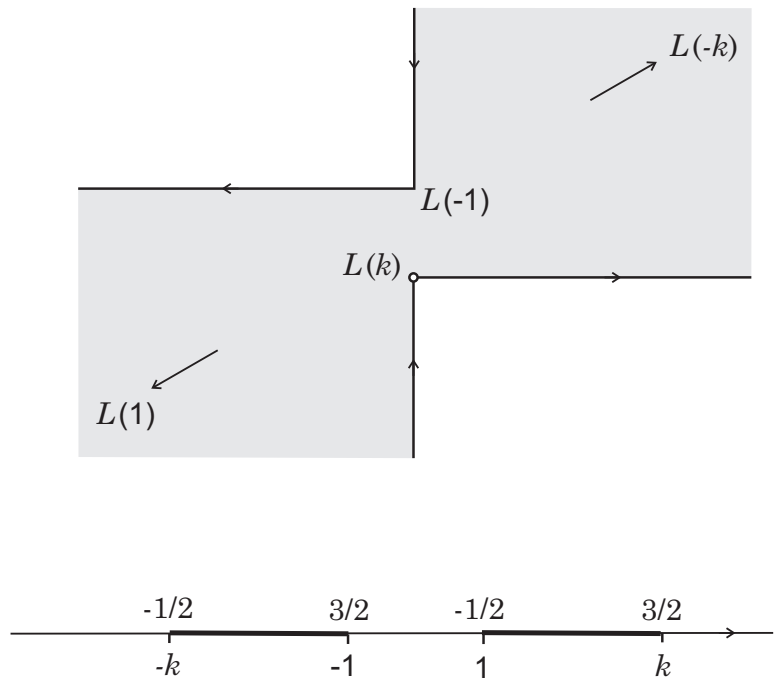}
\caption{Limit of rescaled maps for $k=k_{\mathrm{crit}}$.}
\label{degenerat}
\end{figure}

To find $k_{\mathrm{crit}}$,
we re-scale Fig.~\ref{rect} so that the distance
between $L(-1)$ and $L(k)$ remain fixed.
Consider the rescaled degenerate configuration, Fig.~\ref{degenerat}.
This is a flat quadrilateral.
The map of the upper half-plane on this
region is given by
$$\int\sqrt{\frac{(\zeta-k)(\zeta+1)}{(\zeta+k)(\zeta-1)}}
\frac{d\zeta}{(\zeta+k)(\zeta-1)}$$
which is the same that we obtain if
we put $c=1$ in the integral in $L(z)$.
The condition would be now that
$$\int_{-1}^1\sqrt{\frac{(\zeta-k)(\zeta+1)}{(\zeta+k)(\zeta-1)}}
\frac{d\zeta}{(\zeta+k)(\zeta-1)}=0,$$
but this diverges at one end.
Von Koppenfels \cite[5.3]{vKop}
integrates by parts to obtain his equation for critical $k$.

The modulus of the degenerate region was computed by von Koppenfels
\cite[(5.3.30)]{vKop}
by solving the equation
\begin{equation}\label{K}
K(\kappa')-2E(\kappa')=\int_0^1\frac{dx}{\sqrt{(1-x^2)(1-{\kappa'}^2x^2)}}-2
\int_0^1\sqrt{\frac{1-{\kappa'}^2x^2}{1-x^2}}dx=0.
\end{equation}
The solution is approximately
\begin{equation}\label{kappap}
\kappa^\prime_{\mathrm{crit}}=0.9089085575\ldots.
\end{equation}
The relation between our $k$ and $\kappa'$ is the following
$k=(1+\kappa)/(1-\kappa),$ where $\kappa=\sqrt{1-{\kappa'}^2}.$
This gives the critical value $k_{\mathrm{crit}}=2.4305$.

\vspace{.1in}
Thus as $k$ varies from $1$ to $k_{\mathrm{crit}}$, $c$ varies from $0$ to $1$,
and $K$ (the modulus) varies from $0$ to $K_{\mathrm{crit}}\approx 0.709459$.

When $k\to 1+$, $c\to 0+$; after this point is passed,
we obtain the second family, with modulus $K'$ near infinity.
This second family is symmetric to the first, so there is a correspondence
$K'\mapsto 1/K.$ As seen from the pictures, and using the domain extension
principle we have $K_{\mathrm{crit}}<1$.

Thus for every $K>0$ we have at most one quadrilateral with
angles $(3/2,1/2,3/2,1/2)$. It is of the first type when
$K<K_{\mathrm{crit}}$ and
of the second type when $K>1/K_{\mathrm{crit}}$.
This completes the proof of Theorem 3.

\section{Remarkable constants}

Paper \cite{BE} contains new proofs of some results from \cite{Lin1}.
It is noticed in \cite[Remark 2]{BE} that equation (\ref{K})
and its solution $\kappa^\prime_{\mathrm{crit}}$ in (\ref{kappap})
is a famous constant discussed in Finch \cite[4.5]{Finch}.
To conform to Finch's notation, let
$$c:=\kappa^\prime_{\mathrm{crit}},$$
and
$$\Lambda=\exp\left(\frac{-\pi K\left(\sqrt{1-{c}^2}\right)}{K(c)}\right)=0.1076539192\ldots.$$
There was a constant known for some time as the ``One-Ninth'' constant in
approximation theory, until Gonchar and Rakhmanov proved that it equals
$\Lambda$. Earlier the same constant $\Lambda$
implicitly appeared in Euler's work
on elastics.
Halphen \cite[p. 287]{Halfen} mentions $c$
as the unique solution of the equation
$$\sum_{n=0}^\infty(2n+1)^2(-x)^{n(n+1)},\quad 0<x<1,$$
and computes this solution to six digits.

Related constant
in \cite{Lin1},
$$b_1=\frac{K(c)}{K\left(\sqrt{1-{c}^2}\right)}$$
appears in \cite{Lin1} in the study of metrics of
curvature $1$ with one conic singularity on real tori, but the real tori in
\cite{Lin1} are {\em different} from the real tori studied here.

We recall that there are two types of real tori:
in the first type, which we called rectangular in this paper,
the parameter $\tau=\omega_3/\omega_1$ is pure imaginary.
The second type, which is studied in \cite{Lin1} can be called ``rhombic'',
and can be characterized by the property $\Rea\tau=1/2$.

The problem considered in \cite{Lin1} has no solutions for rectangular tori,
while it has a unique solution for the rhombic tori with parameter
$\tau=1/2+ib$ exactly when $b>b_1$ or $b<1/(4b_1)$.

That the same constant arises in these two different settings suggest some
relation between the problem studied in \cite{Lin1} and the
problem studied here in Section 4\footnote{In a private communication,
C.-S. Lin explained this coincidence: there is a degree $2$ isogeny
$T\to T'$ between a rectangular torus $T$ with primitive periods $(1,it)$ and
rhombic torus $T'$ with primitive periods $(1\pm it)/2$. Every metric
on $T'$ with one conic singularity at $0$ pulls back via this isogeny
to a metric on $T$ with singularities at $0$ and $(1+it)/2$. Using this
correspondence, Proposition 1 can be derived from the results in \cite{Lin1},
but our proof of Proposition 1 is much more elementary.}

\section{Algebraic developing maps}

When the angle $\alpha$ in Fig.~\ref{2families} is rational, the monodromy
of $f$ is a finite group, therefore $f$ is an algebraic function.
In this case,
the general solution of (\ref{heun}) is algebraic, and according to a theorem
of Klein, (\ref{heun}) is a pull-back of a hypergeometric equation
via a rational function \cite{BD}. This type of Heun's
equations was intensively
studied, see, for example \cite{FV1,FV2,V}, but the case of algebraically
solvable Heun's equation considered here seems to be new.

The procedure of obtaining an explicit expression of $f$
is explained in detail in \cite{EG3}. Let $\alpha=p/q,\; (p,q)=1$. 
The monodromy group of the developing map $f$
is a dihedral group of order $2q$. It is generated by
$z\mapsto 1/z$ and $z\mapsto e^{2\pi i\alpha}z$.
The fundamental invariant of this dihedral group is
the rational function 
$$g(z)=-\frac{1}{4}\left( z^q+\frac{1}{z^q}-2\right).$$
This is a Belyi function which means that its critical values
are $0,1,\infty$. It is real on the real line, on the unit circle
and on the line $\ell_\alpha$. It is easy to see that $h=g\circ f$
is also a Belyi function, and for small values of $q$ function
$h$ can be found explicitly, which gives an explicit expression of $f$.

Here we only give few examples.
\vspace{.1in}

\noindent
Example 1.  $\alpha=1/2$. In this case, the image
of the real line under $f$ is contained in three perpendicular circles:
the real and imaginary lines and the unit circle.

Consider the function
$$g(z)=-\frac{1}{4}\left(z-\frac{1}{z}\right)^2.$$
The critical points $0,\infty,\pm1,\pm i$ are all simple, and the critical
values are
$$g(0)=g(\infty)=\infty,\quad g(\pm1)=0,\quad g(\pm i)=1.$$
Function has the symmetry properties
$$g(\overline{z})=g(-\overline{z})=g(1/\overline{z})=\overline{g(z)}.$$
Therefore, real and imaginary lines and the unit circle are mapped on
the real line, and these three circles constitute the full preimage
of the real line. These three circles define a cell decomposition of
the sphere into $8$ triangles, each is mapped bijectively onto the upper
or lower half-plane. Function $g$ is totally ramified over $0,1,\infty$,
with simple ramification points over each.

Consider the function $h=g\circ f$ in the upper half-plane.
It
is real on the real line, so it extends by symmetry to a rational
function. This rational function has the property that it is ramified
only over $0,1,\infty$, so it is a Belyi function \cite{Schneps}.
It is not difficult to write $h$ explicitly:
$$h(z)=-\frac{(z+2)(z-2)^3}{3(z^2+2z-2)^2}.$$
The critical points are:

$-1$ of multiplicity $2$ with critical value $1$,

$2$ of multiplicity $2$ with critical value $0$,

and two simple critical points $-1\pm\sqrt{3}$ with critical value $\infty$.

It follows that $f=g^{-1}\circ h$ maps the upper half-plane onto
a rectangle with angles $(1/2,3/2,1/2,3/2)$, as in Fig.~\ref{2families}
with $\alpha=1/2$, and $k=2$, $K=0.63963$
in this case.
\vspace{.1in}

\noindent
Example 2. $\alpha=1/3$. Monodromy group is dihedral $D_3$ in this case, and we can
find $f$ in the form $f=g^{-1}h$ where
\begin{equation}\label{1/3}
g(z)=-\frac{1}{4}\left(z^3+\frac{1}{z^3}-2\right),
\end{equation}
and $h$ is a Belyi function of the form
$$h(z)=s\frac{(z-x)^2(z-a)}{(z-y)^3(z-t)^3},$$
with parameters $a,x,y,s,t$ chosen so that the following conditions
are satisfied:
$$h(0)=h(1)=h(w)=1,\quad h'(1)=h''(1)=h'(w)=0.$$
Computation with Maple gives the following exact values, where
we denoted $\epsilon=2^{1/3}>0$:
$$a=\frac{5}{4}\epsilon^2+\frac{3}{2}\epsilon+3\approx 6.87,$$
$$x= -\frac{1}{10}\epsilon^2-\frac{3}{10}\epsilon+\frac{3}{5},$$
$$y=\frac{1}{2}\epsilon^2+\frac{1}{2}\epsilon+\frac{3}{2}-\frac{1}{2}
\sqrt{8\epsilon^2+10\epsilon+13},$$
$$t=\frac{1}{2}\epsilon^2+\frac{1}{2}\epsilon+\frac{3}{2}+
\sqrt{8\epsilon^2+10\epsilon+13},$$
$$s=\frac{33}{4}\epsilon^2+\frac{21}{2}\epsilon+13.$$
$$w=\frac{3}{2}\epsilon^2+\frac{3}{2}\epsilon+3.$$
Function $h$ has critical values $0,1,\infty$ whose preimages are:
\vspace{.1in}

Zeros: $a$ of multiplicity $1$, $x$ of multiplicity $2$ and
$\infty$ of multiplicity $3$.

$1$-points: $0$ of multiplicity $1$, $w$ of multiplicity $2$ and
$1$ of multiplicity $3$.

Poles: $y$ and $t$ both of multiplicity three.

By a direct verification $g^{-1}(h)$ maps the upper half-plane
onto a rectangle equivalent to that shown
in Fig.~\ref{2families} with $\alpha=1/3.$

We have $K\approx 0.67957$ for this quadrilateral.
\vspace{.1in}

\noindent
Example 3. $\alpha=2/3$. Function $g$ is the same as in (\ref{1/3}), and
consider the Bely function
$$h=\frac{64(135+78\sqrt{3})(z-1)^3}{(z-4-2\sqrt{3})^3(3z+2\sqrt{3})^3}.$$
Equation $h(z)=1$ has two simple real solutions, $0$ and $8+4\sqrt{3}$,
and two complex conjugate solutions of multiplicity $2$.
Function $h$ has $2$ triple zeros, at $1$ and $\infty$,
and two triple poles. Then $f=g^{-1}\circ h$ is the developing map
of a spherical rectangle equivalent to that shown in
Fig.~\ref{2families}.
Parameter corresponding to this example is  $K=0.57735$.

{\em Department of Mathematics, Purdue University,
West Lafayette Indiana, 47907 USA}
\end{document}